\documentclass[12pt,a4paper]{article}

\pdfoutput=1

\usepackage{indentfirst}
\usepackage{enumerate}
\usepackage[pdftex]{graphicx}
\usepackage{graphicx}
\usepackage{epstopdf}
\usepackage{amsmath}
\usepackage{amssymb,amsfonts}
\usepackage{multicol}
\usepackage{url}
\usepackage{color}
\definecolor{DarkGreen}{rgb}{0,0.5,0}
\definecolor{DarkRed}{rgb}{0.5,0, 0}

\newcommand{\5}{\mbox{d\llap{\raise0.75ex\hbox{-}}}}
\newcommand{\6}{\mbox{\rlap{\raise0.3ex\hbox{--}}D}}
\newcommand{\4}{\mbox{\rlap{\raise0.15ex\hbox{-}}d}}

\newtheorem{theorem}{Theorem}
\newtheorem{definition}{Definition}
\newtheorem{lemma}{Lemma}

\newtheorem{proposition}{Proposition}

\begin{document}

\centerline {\Large\bf $q-$spherical surfaces in Euclidean space}

\bigskip

 \begin{center}

 SONJA GORJANC\\
{\small University of Zagreb, Faculty of Civil Engineering,\\ Ka\v{c}i\'{c}eva 26, 10000 Zagreb, Croatia\\
e-mail: sgorjanc@grad.hr}

\medskip

 EMA JURKIN\\
{\small University of Zagreb, Faculty of Mining, Geology and Petroleum Engineering,\\ Pierottijeva 6, 10000 Zagreb, Croatia\\
e-mail: ema.jurkin@rgn.hr}
\end{center}

\bigskip

\noindent {\small {\bf Abstract}.\\
In this paper we define $q$-spherical surfaces as the surfaces that contain the absolute conic of the Euclidean space as a $q-$fold curve. Particular attention is paid to the surfaces with  singular points of the highest order. Two classes of such surfaces,  with one and two $n-$fold  points, are discussed in detail. We study their properties, give their algebraic equations and visualize them with the program {\it Mathematica}.

\medskip
\medskip

}

\bigskip

\section{Motivation}

One of the most important places in the classical geometry belongs to the study of some classes of surfaces with special properties in the Euclidean space.
 In the real projective plane the Euclidean metric defines the Euclidean plane $\mathbb E^2$  with the \textit{absolute}  (\textit{circular}) \textit{points} $(0,1,i)$ and $(0,1,-i)$ (i.e.  $x_0=0$, $x_1^2+x_2^2=0$). An algebraic curve of order $n$ passing through the absolute points is called \textit{circular curve}. If it contains absolute points as its $q-$fold points, the curve is called $q-$\textit{circular}, and if $n=2q$ the curve is called \textit{entirely circular}. Due to their numerous applications in engineering, the circular curves in the Euclidean plane (like strophoid, trisectrix of Maclaurin, lima\c con, cardiod, lemniscate of Bernoulli, Booth lemiscate, Cartesian ovals, Cassini ovals, astroid, Watt's curve) have been treated in a significant number of papers (see e.g. \cite{basset}, \cite{fladt62}, \cite{savelov}, \cite{wieleitner}, \cite{wiki-circular}, \cite{mathworld}). On the other hand, the surfaces in the Euclidean space with the similar properties (that contain the absolute conic of $\mathbb E^3$ as a $q-$fold curve), with an exception of the cyclides that will be mentioned later, have been in some way neglected. Motivated by this fact, we introduce so-called $q$-spherical surfaces, study some of their properties and visualize their shapes. 
Particular attention will be paid to the surfaces with  singular points of the highest order.
  
\section{Introduction}

In the real three-dimensional projective space $P^3(\mathbb{R})$, in homogeneous Cartesian coordinates $(x_0,x_1,x_2,x_3)$, ($x_1,x_2,x_3\in \mathbb{R}, x_0\in\{0,1\}, (x_0,x_1,x_2,x_3) \neq (0,0,0,0)$), the equation 
$$F_{n}(x_0,x_1,x_2,x_3)=0,$$ 
where $F_{n}$ is a homogeneous  algebraic polynomial of degree $n$, defines an $n^{th}$ order surface  $S_{n}$. This equation can be written as 
\begin{equation}\label{surface}
f_n(x_1,x_2,x_3)+x_0f_{n-1}(x_1,x_2,x_3)+ ... +x_0^{n-1}f_1(x_1,x_2,x_3)+x_0^{n}f_0(x_1,x_2,x_3)=0,
\end{equation}
 where $f_{j}$, $j=1, ..., n$, are homogeneous  algebraic polynomials of degree $j$, \cite{Salmon1}. 

Any straight line,  not lying on  $S_{n}$, intersects  $S_{n}$ in $n$ points and any plane intersects  $S_{n}$ in the  $n^{th}$ order plane curve. 

A point $T$ of the surface $S_{n}$ for which at least one partial derivation of $F_{n}$ is not equal to zero is called the regular point of  $S_{n}$. All tangents to the surface at that point lie in one plane - the tangent plane of  $S_{n}$ at $T$. 

A point $T$ of the surface $S_{n}$ for which all partial derivations of $F_{n}$ are equal to zero is called the singular point of  $S_{n}$. The tangents to  $S_{n}$  at this point form an algebraic cone with vertex in $T$. If the tangent cone is of order $k$, the point $T$ is the $k$-fold point of the surface   $S_{n}$. Every  plane through $T$  intersects  $S_{n}$ in the  $n^{th}$ order plane curve with the $k$-fold point in $T$. 

The point $T$ is a $k-$fold  point of $S_n$ if all partial derivations of $F_n$ with the order less than $k$ vanish at $T$, and at least one $k-$order derivative of $F_n$ at $T$ doesn't vanish, \cite{Goetz}.

According to \cite{Harris}, if the origin $O(1,0,0,0)$ is the $k$-fold point of   $S_{n}$, then  $S_{n}$ has the equation,  
\begin{equation}\label{k-fold}
f_n(x_1,x_2,x_3)+x_0f_{n-1}(x_1,x_2,x_3)+ ... +x_0^{n-k}f_k(x_1,x_2,x_3)=0,
\end{equation}
and the tangent cone at $O$ is given by
\begin{equation}
f_k(x_1,x_2,x_3)=0.
\end{equation}
\medskip

In the real projective space $P^3(\mathbb{R})$ the Euclidean metric defines the Euclidean space $\mathbb E^3$  with the absolute conic given by the equations: $$ x_0=0\,\, {\textnormal {and}}\,\,A_2=x_1^2+x_2^2+x_3^2=0.$$

\section{$q-$spherical surfaces}

\begin{definition}
A surface $S_n$ of Euclidean space is called $q-$spherical surface if it contains the absolute conic as a $q-$fold curve.
\end{definition}

\begin{theorem} In Euclidean space $q-$spherical surface of order $n$ is given by the following equation:
\begin{equation}\label{formula}
A_2^q\, g_{n-2q}(x_1,x_2,x_3)+\sum_{j=1}^{q-1}x_0^j A_2^{q-j} g_{n-2q+j}(x_1,x_2,x_3)+\sum_{j=q}^{n}x_0^j  f_{n-j}(x_1,x_2,x_3)=0,
\end{equation}
where $n\ge 2q$,  $g_{n-2q}\neq 0$, $A_2\nmid g_{n-2q}$, $f_{n-q}\neq 0$ and $A_2\nmid f_{n-q}$.
\end{theorem}
{\sc Proof:}
Let us first prove that the surface given by (\ref{formula}) is  $q-$spherical surface.   For all points on the absolute conic ($A_2=0,\, x_0=0$) all first $q-1$ partial derivatives of the polynomial on the left hand side of (\ref{formula})  vanish since their terms contain either  $A_2$ or $x_0$ as a factor. According the theorem's condition, at least one $q-$order derivative does not vanish for the points of the absolute conic. 

Now, let a $q-$spherical surface of order $n$ be given by the equation (\ref{surface}) that can be written in the following form:
\begin{equation}
F_n(x_0,x_1,x_2,x_3)=f_n(x_1,x_2,x_3)+\sum_{j=1}^{q-1}x_0^j f_{n-j}(x_1,x_2,x_3)+\sum_{j=q}^{n}x_0^j  f_{n-j}(x_1,x_2,x_3)=0.
\end{equation}

Since the first $q-1$ partial derivations of the polynomial  $F_n$ have to vanish for $A_2=0$ and $x_0=0$, it is clear that the polynomials $f_{n-j}$, $j=1,\dots ,q-1$, must contain $A_2^{q-j}$ as a factor. Since at least one $q-$order derivative of $F_n$ doesn't vanish, it is clear that $f_{n-q}\neq 0$ and $A_2$ is not a factor of $f_{n-q}$. 
\hfill{$\square$}  

\medskip

In  \cite{moniods-nasi} and \cite{quartics}  the authors studied the special class of 1-spherical surfaces, the surfaces  which touch the plane at infinity through the absolute conic and have a singular point of the highest order. These surfaces belong to a wider class of surfaces, so-called monoid surfaces, also treated in \cite{monoids} and introduced  in \cite{rohn}.

\bigskip

In  paper \cite{circular} the authors considered a congruence of circles $\mathcal C(p)$ that consists of circles in Euclidean space $\mathbb E^3$  passing through two given points $P_{1,2}(0,0,\pm p)$. For a given congruence $\mathcal C(p)$ and a given curve $\alpha$,  a {\it circular surface} $\mathcal {CS}(\alpha,p)$ is defined as the system of circles from $\mathcal C(p)$  that intersect $\alpha$. If $\alpha$ is an $m^{th}$ order algebraic curve that cuts the axis $z$ at $z'$ points, the absolute conic at $a'$ pairs of the absolute points and with the points $P_1$ and $P_2$ as $p'_1$-fold and $p'_2$-fold points, respectively, then, the following statements hold:\\
\hspace*{0.5cm} -- $\mathcal{CS}(\alpha,p)$ is an algebraic surface of the order $3m-(z'+2a'+2p'_1+2p'_2)$.\\
\hspace*{0.5cm} -- The  absolute conic is an  $m-(z'+p'_1+p'_2)$-fold curve of $\mathcal{CS}(\alpha,p)$.\\
\hspace*{0.5cm} -- The axis $z$ is an $(m-2a'+z')$--fold line of $\mathcal{CS}(\alpha,p)$.\\
\hspace*{0.5cm} -- The points  $P_1$, $P_2$ are  $2m-(2a'+p_1'+p_2')$--fold points of \,$\mathcal{CS}(\alpha,p)$.\hfill ({\bf CS})
\\
Thus, the most of these surfaces are spherical, and some of them will be studied in more detail in  subsection \ref{potpoglavlje}. 

\bigskip

The further examples of $q-$spherical surfaces can be found in \cite{roses}  and \cite{generalized_roses} where authors introduced  rose surfaces and generalized rose surfaces as the special cases of circular  surfaces for which $\alpha$ is a cyclic-harmonic curve. It was shown how $q$ depends on the properties and  position of the referent cyclic-harmonic curve with respect to the singular points of the congruence  $\mathcal C(p)$.

\section{Entirely spherical surfaces}

\begin{definition}
A surface $S_{2n}$ of Euclidean space is called entirely spherical surface if it contains the absolute conic as an $n-$fold curve.
\end{definition}

Probably the mostly studied surfaces of the fourth order are cyclides, the bispherical quartic surfaces, \cite{Salmon2}.
The term ``cyclides'' is often used for their special class, Dupin cyclides, which can be defined in a several ways, \cite{eisenhart}, \cite{hilbert}. They are only surfaces that have the property that their evolulute  degenerate into curves,  in fact both sheets of the focal surface are conics. The examples of cyclides are the tori, cones and cylinders of revolutions. 
   Dupin cyclides can also be defined as the envelopes of the family of spheres tangent to three fixed spheres.  
 Dupin cyclides are the inverse images of the standard tori, cylinders or cones. The inverse image of a ring torus, horn torus and spindle torus are called a ring cyclide, horn cyclide and spindle cyclide, respectively. If the center of the inversion sphere lies on the torus, the obtained surface is a parabolic cyclide (ring, horn or spindle), \cite{math-cyc}. 
The cyclides as the surfaces of the forth order having the circle at infinity as the nodal conic were studied in \cite{Salmon2}.  Their equation of type (\ref{formula}) for $n=4, p=2$, was given. The cyclide was also defined as the envelope of sphere whose center moves on a fixed quadric, and which cuts a fixed sphere orthogonally. The intersection curve of the fixed quadric and sphere is the focal curve of the cyclide. The number of strait lines lying on the cyclide is sixteen. In \cite{Salmon2} the author distinguishes 23 types of cyclides.
In \cite{pottmann} the family of surfaces called Darboux cyclides were studied. These surfaces are algebraic surfaces of order at most 4 and are a superset of bispherical surfaces of order 4, circular surfaces of order 3 and quadrics, and they carry up to six real families of circles. 

\medskip

Entirely spherical surfaces can be constructed in many different ways. Let us mention some of them: 

An inverse image of a surface $S_n$ of order $n$ is an $n-$spherical surface $S_{2n}$ of order $2n$ with an $n-$fold point in the pole of inversion.

If $C_n$ is a surface of class $n$ and $P$ a point in general position to $C_n$, then the pedal surface of $C_n$ with respect to the pole $P$ is an $n-$spherical surface $S_{2n}$ of order $2n$ with an $n-$fold point in the pole $P$, \cite{kranj}.

According to the properties ({\bf CS}) of circular surfaces $\mathcal {CS}(\alpha,p)$ that are pointed out in the previous section,  if $\alpha=k_{2n}$ is an entirely circular curve of order $2n$ with an $n-$fold point in $P_1$ ($m=2n$, $a'=p_1=n$, $p_2=0$, $z'=0$),  the obtained surface $\mathcal {CS}(k_{2n},p)$ is entirely spherical surface with two $n-$fold points $P_1, P_2$. The examples of these surfaces will be given in  subsection \ref{potpoglavlje}.

\begin{theorem}
An entirely spherical surface $S_{2n}$ of order $2n$ can't have singular points of multiplicity higher then $n$.
\end{theorem}
{\sc Proof:} If there was a point $T$ of multiplicity  $n+1$, all isotropic lines through $T$ would lie on $S_{2n}$ and the surface would split onto the isotropic cone with vertex in $T$ and a surface of order $2n-2$. 
\hfill{$\square$}

\begin{lemma} \label{lemma1}
A plane curve of order $2n$, $n\geq 2$, can have at most three $n-$fold points.
\end{lemma}
{\sc Proof:} The maximum number of double points of a curve of order $2n$ equals $D=\frac{(2n-1)(2n-2)}{2}$, \cite{SalmonCurves}.  Every $n-$fold point is counted as $N=\frac{n(n-1)}{2}$ double points. Since $\frac{D}{N}=4-\frac{2}{n}$, the curve of order $2n$ can have three, but not four $n-$fold points.   
\hfill{$\square$}  

\medskip

Therefore, if a plane curve  of order $2n$  has more then three $n-$fold points, it splits onto the curves of lower order.   
 
\begin{lemma} \label{one_nfold}
An entirely circular curve of order $2n$, $n\geq 2$, can have only one $n-$fold point beside the absolute points. 
\end{lemma}
{\sc Proof:} The statement follows directly from Lemma \ref{lemma1}.
\hfill{$\square$}

\begin{theorem}\label{circ_surf}
If an entirely spherical surface $S_{2n}$ contains two real and different $n-$fold points $N_1$ and $N_2$, every plane
through the line $N_1N_2$ cuts $S_{2n}$ into $n$ circles passing through $N_1$ and $N_2$.
\end{theorem}
{\sc Proof:} Let  $N_1$ and $N_2$ be the $n-$fold points of the surface $S_{2n}$, and let $\beta$ be a plane through $N_1$ and $N_2$. The intersection line of $\beta$ and $S_{2n}$ is a curve $k_{2n}$ of order $2n$. According to Lemma \ref{one_nfold} this curve is not proper, it splits onto curves $k_{2t}$ and $k_{2n-2t}$ of order $2t$ and $2n-2t$, respectively. Since $k_{2n}$ passes through the absolute points $n$ times, the curves $k_{2t}$, $k_{2n-2t}$ have to be entirely circular passing through the absolute points $t$ and $n-t$ times. The necessary multiplicity of the point $N_1$ will be achieved if it is the singular point of $k_{2t}$ and $k_{2n-2t}$ with multiplicity $t$ and $n-t$. Now, again according to Lemma \ref{one_nfold}, $N_2$ is the singular point of $k_{2t}$ and $k_{2n-2t}$ with multiplicity at most $t-1$ and $n-t-1$. This is in contradiction with the fact that  $k_{2n}$ passes  $n$ times through $N_2$. It follows that $k_{2t}$ and $k_{2n-2t}$ are not proper curves either. Therefore, they again split onto the curves of lower order. Continuing with this procedure, in the last step we come to the $n$ curves of order 2, i.e. circles.      
       \hfill{$\square$}

\subsection{Examples}

In this subsection we give some examples of the spherical surfaces. For a particular  surface given by an implicit equation, we study its properties and visualize its shape. For  computing and plotting, we use  the program {\it Mathematica}.
In some examples we start with an algebraic equation of the surface and determine its properties. In other cases we first give a construction of the surface from which we derive its equation.

\subsubsection{Entirely spherical surfaces $S_{2n}$ with only one $n$-fold point}
A surface $S_{2n}$ given by the equation of the form 
\begin{equation}
A_2(x_1,x_2,x_3)^n+x_0^n\cdot f_{n}(x_1,x_2,x_3)=0, \nonumber
\end{equation}
i.e. in the affine coordinates
\begin{equation}\label{one_nfold_surf}
A_2(x,y,z)^n+f_{n}(x,y,z)=0,
\end{equation}
is an entirely spherical surface with the $n-$fold point at the origin.
\begin{proposition}
There are only $2n$ straight lines through the origin lying entirely on the surface $S_{2n}$  given by (\ref{one_nfold_surf}). These lines are the intersections of cones given by
\begin{equation}
 A_2(x,y,z)=0, \quad f_n(x,y,z)=0,   \nonumber
\end{equation}
and they are imaginary in pairs.
\end{proposition}
{\sc Proof:} 
Let a line $p$ through $O(0, 0, 0)$ be spanned by $O$ and a further point $P(a,b,c) \neq O$.  The line $p$ is parametrized by
\begin{equation}\label{pravacp}
p \quad ... \quad (x,y,z)=(ta, tb, tc), t \in \mathbb{R}.  
\end{equation}
It lies on $\mathcal {S}_{2n}$ if and only if  \vspace*{-0.1cm}
\begin{equation}
A_2(ta, tb, tc)^n+f_{n}(ta, tb, tc)=0,   \nonumber
\end{equation} \vspace*{-0.1cm}
for every $ t \in \mathbb{R}$.
This is precisely when \vspace*{-0.1cm}
\begin{equation}
t^{n}[t^{n}A_2(a,b,c)^n+ f_{n}(a,b,c)]=0,   \nonumber
\end{equation} 
for every $ t \in \mathbb{R}$.
It follows that 
$A_2(a,b,c)=0$ and   $f_{n}(a,b,c)=0. $
Therefore,  
$A_2(ta, tb, tc)=0$, $f_{n}(ta, tb, tc)=0$, 
for every $ t \in \mathbb{R}$.
Evidently the line $p$ lies on the cones given by equations $A_2(x,y,z)=0$ and $f_n(x,y,z)=0$. We conclude: the only lines through the origin that lie on  $\mathcal {S}_{2n}$ are the isotropic lines on the tangent cone at  the origin.\hfill {\small $\square$}

\begin{proposition}
The surface $S_{2n}$ given by  (\ref{one_nfold_surf}) has no other singular points beside the origin.
\end{proposition}
{\sc Proof:} This fact can be proved as follows: Let  a non-isotropic line $p$ through the origin be given by (\ref{pravacp}). We compute the intersections of $S_{2n}$ and $p$. They belong to the zeros of the following polynomial of degree $2n$ in $t$: 
$$P(t):=A_2(ta, tb, tc)^n+ f_{n}(ta, tb, tc).$$
Obviously, $t=0$ is zero with multiplicity $n$. The other $n$ zeros are given by $$t^n=-\frac{f_{n}(a, b, c)}{(a^2+b^2+c^2)^n}$$ and therefore different (in general) complex numbers.\hfill {\small $\square$}
  
\bigskip

The tangent cone $\Phi_n$ at the $n-$fold point $O$ of the surface $S_{2n}$ given by (\ref{one_nfold_surf}) has the equation  $f_{n}(x,y,z)=0$. The polynomial $f_{n}$ can be irreducible or reducible. If the polynomial $f_{n}$ can be factorized
\begin{equation}\label{one_nfold_fact}
f_{n}(x,y,z)=f_{n_1}(x,y,z)\cdot ... \cdot f_{n_k}(x,y,z), \quad n_1+...+n_k=n,
\end{equation}
the tangent cone $\Phi_n$ splits into the cones $\Phi_{n_1}, ..., \Phi_{n_k}$ of order $n_1, ..., n_k$, respectively. Therefore, the classification of surfaces $S_{2n}$ can be made according to the degrees of the polynomials $f_{n_1}, ...,  f_{n_k}$. 
If we assume that all polynomials $f_{n_1}, ...,  f_{n_k}$ determine different cones $\Phi_{n_1}, ..., \Phi_{n_k}$, then for all $n\in N$ the surfaces $S_{2n}$ can be classified into $p(n)$ types, where $p$ is the partition function, i.e. $p(n)$ is the number of ways of writing the integer $n$ as a sum of positive integers, where the order of addends is not considered significant. The formulas for counting the value $p(n)$ can be found in \cite{math-partition}.

Some examples of surfaces given by (\ref{one_nfold_surf}), together with their tangent cones at the $n-$fold points, are shown in the next four figures.

\begin{center}
 \includegraphics[width=15cm]{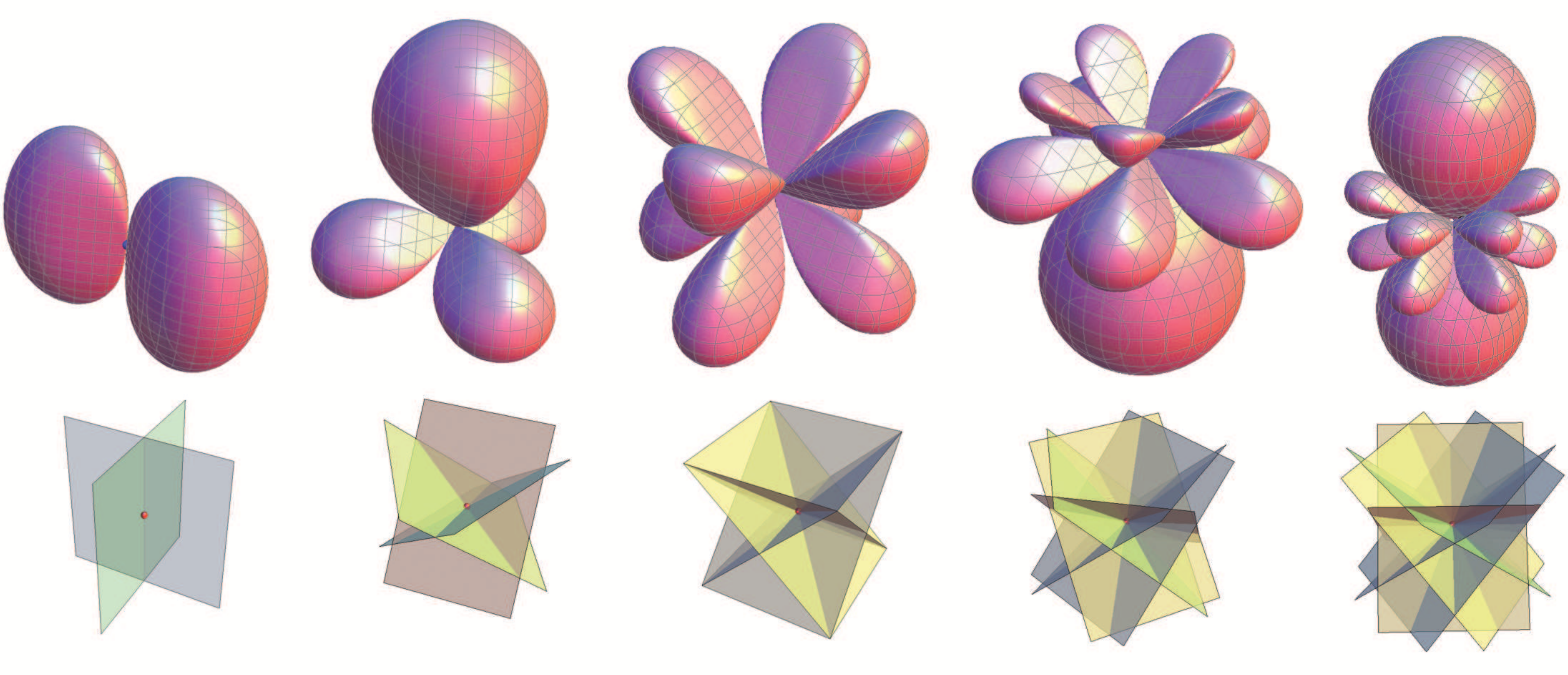}\\
\vspace{-0.5cm}
\begin{multicols}{5}
a

b

c

d

e

\end{multicols}
\end{center}
Figure 1: Five examples  for $n=2,3,\dots ,6$ where the tangent cone at the origin splits into $n$ planes, i.e. $f_n(x,y,z)=\prod_{i=1}^n h_i(x,y,z) $, where $\forall i$,  $h_i(x,y,z)$ is a linear polynomial.

\begin{center}
 \includegraphics[width=15cm]{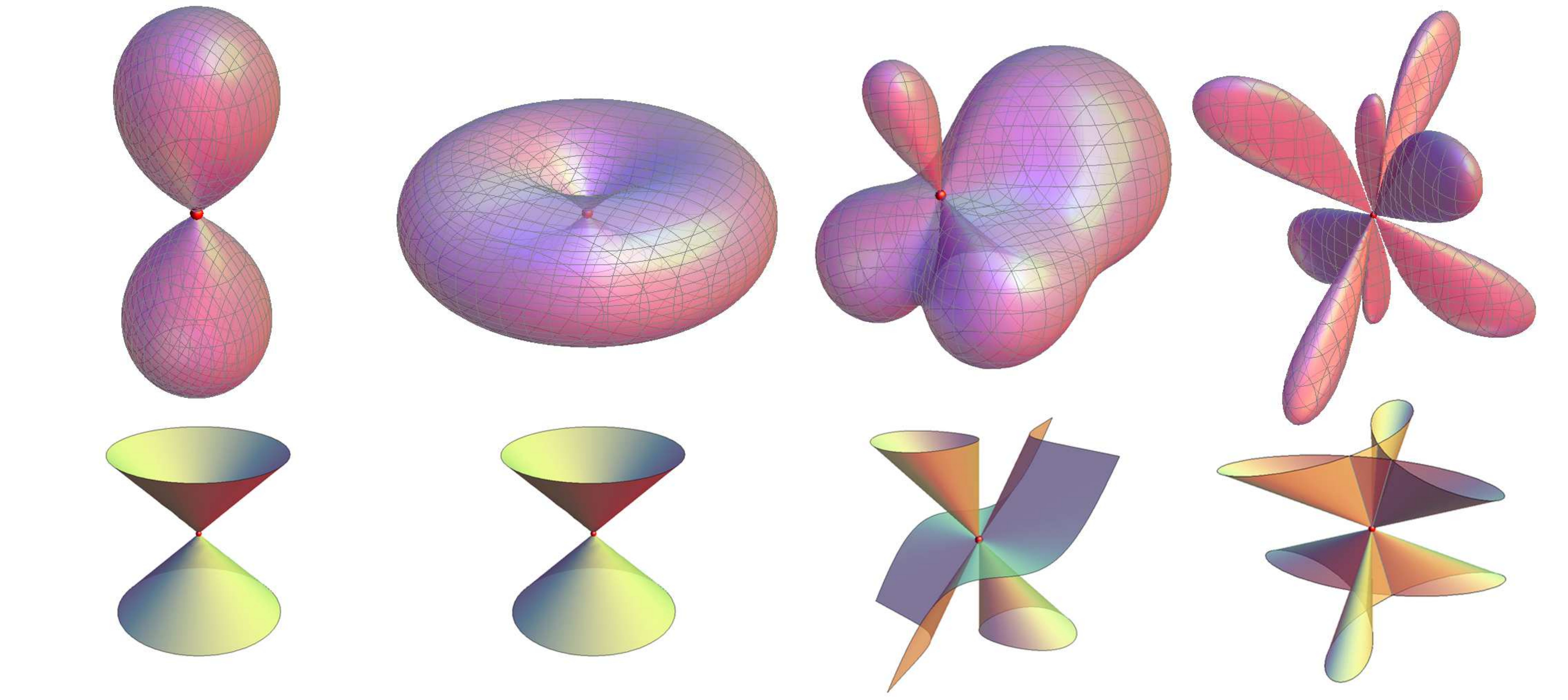}
\\
\begin{multicols}{4}
a

b

c

d

\end{multicols}
\end{center}
Figure 2: Surfaces $S_{2n}$ with the  proper  tangent cones of degree $n$ at the origin, where  $f_n(x,y,z)$ is:
$x^2 + y^2 - z^2$ (a),
 $-(x^2 + y^2 - z^2)$ (b),
 $-2 x^3 - x^2 z + 2 y^2 z + x z^2$ (c),
 $48 x^4 + 48 y^4 - 64\sqrt{3} y^3 z + 40 y^2 z^2 - z^4 + 
 8 x^2 (12 y^2 + 24 \sqrt{3} y z + 5 z^2)$ (d).

\begin{center}
 \includegraphics[width=15cm]{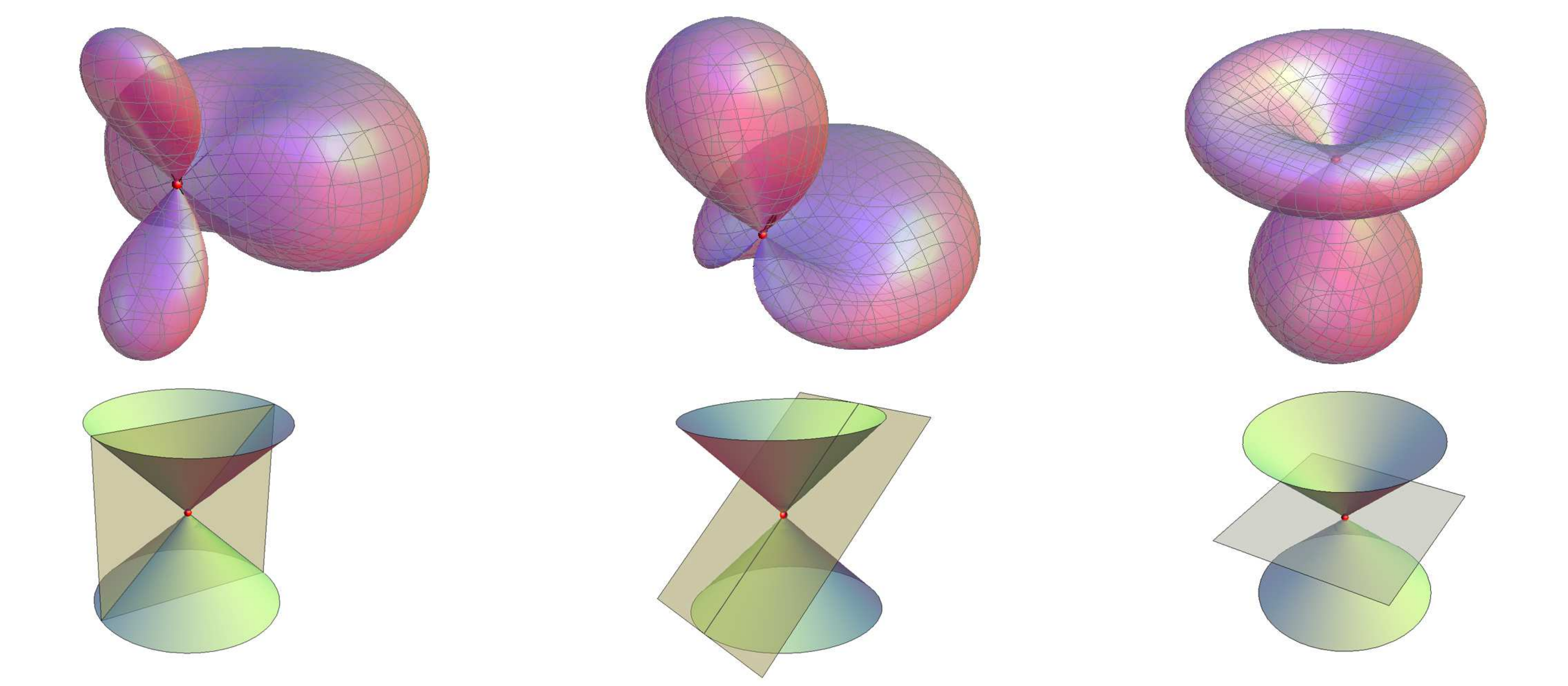}\\
\vspace{-0.5cm}
\begin{multicols}{3}
a

b

c

\end{multicols}
\end{center}
Figure 3: Three examples  for $n=3$ where the tangent cone at the origin splits into a plane and 2nd degree cone.

\begin{center}
 \includegraphics[width=15cm]{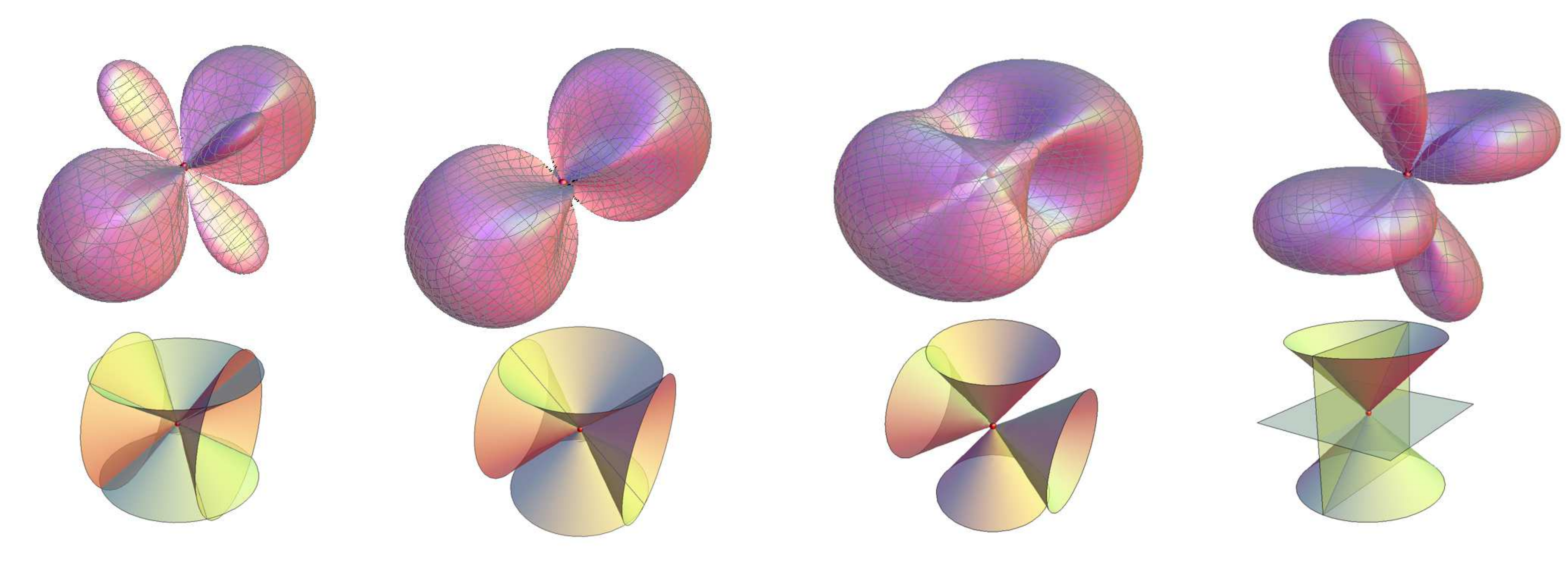}\\
\vspace{-0.5cm}
\begin{multicols}{4}
a

b

c

d

\end{multicols}
\end{center}
Figure 4: Four examples  for $n=4$ where the tangent cone at the origin splits into two 2nd degree cone or into two planes and one 2nd degree cone.

\subsubsection{Entirely spherical surfaces $S_{2n}$ with two $n$-fold points}\label{potpoglavlje}

In this subsection we study a class of entirely spherical surfaces having two singular points of the highest order. We start with defining a class of entirely circular curves.  

\smallskip

A curve  $k^{2n}$ of order $2n$, $n \geq 2$, given by the equation 
\begin{equation}\label{circular curve}
\left(x^2+y^2 \right)^n + f_{n}(x,y) =0, 
\end{equation}
where homogeneous  algebraic polynomial $f_{n}$ is a product of $n$ linear factors, is an entirely circular  curve having an $n-$fold point at the origin. Linear factors of $f_n$ represent the tangent lines at the singular point. If the tangent lines divide the plane into equal parts, polynomial $f_{n}$ equals 
$$
f_n= \left\{\begin{array}{ll}
{\displaystyle\prod_{i=0}^{n-1} \left(\cos i\frac{2\pi}{n} \cdot y - \sin i\frac{2\pi}{n} \cdot x       \right)}, & n \quad \textrm{odd},\\
{\displaystyle\prod_{i=0}^{n-1} \left(\cos i\frac{\pi}{n} \cdot y - \sin i\frac{\pi}{n} \cdot x       \right)}, & n \quad \textrm{even}.
\end{array} \right.
$$

Some examples of this type of curves are shown in Figure~5.

\begin{center}
 \includegraphics[width=15cm]{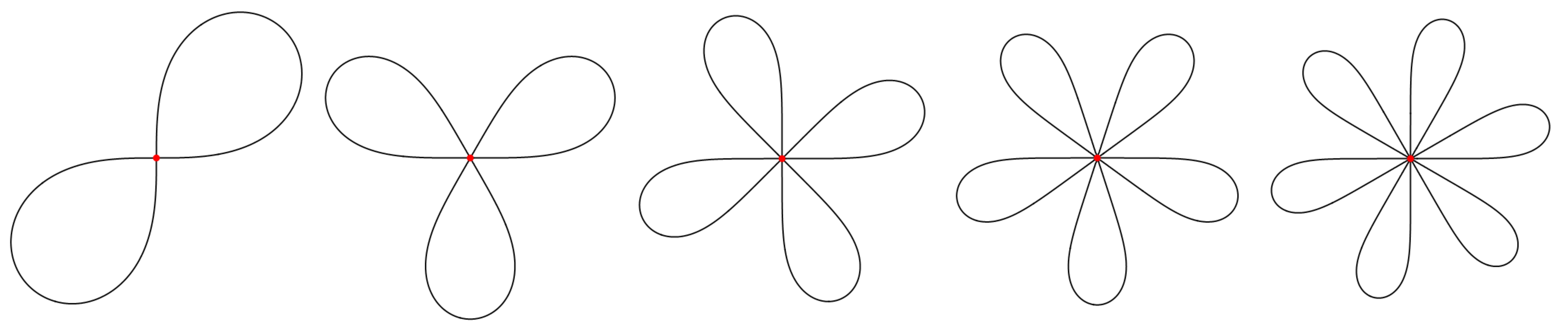}\\
\vspace{-0.5cm}
\begin{multicols}{5}
a

b

c

d

e

\end{multicols}
\end{center}
Figure 5: Five examples  of entirely circular curves with equation \eqref{circular curve} for $n=2, ..., 6$.

\medskip
{\bf Remark:}
The examples in Figure 5 make it easy to see that every straight line through the origin intersects  $k^{2n}$, except at the $n-$fold point $O$, at just one real point  if $n$ is odd number, and that the number of real intersections is zero or two if $n$ is even.  This fact can be also numerically verify for any chosen $n$ and straight line through $O$ using the program {\it Mathematica}.

Switching to polar coordinates  $O(\rho,\varphi)$, where $x=\rho\cos\varphi$ and $y=\rho\sin\varphi$, we obtain the following polar equation of the curve $k^{2n}$

\begin{equation}\label{polar - curve}
\rho^n= (-1)^n 2^{1-n} \sin (n \varphi), \quad \varphi\in [0,2\pi).
\end{equation}

Without lost of generality, from now on we  assume that 
polar equation of $k^{2n}$ is given by expression:

\begin{equation}\label{polar - curve1}
\rho=\sqrt[n]{\sin(n\varphi)}, \quad \varphi\in [0,2\pi).
\end{equation}

\bigskip

\bigskip
Let us now considered a congruence  $\mathcal C(p)$ that consists of circles passing through two given points $P_{1,2}(0,0,\pm p)$, where $p$ is a positive real number. 
For the congruence $\mathcal C(p)$ and  curve $k^{2n}$ given by \eqref{circular curve},  a  circular surface $\mathcal {CS}(k^{2n}, p)$ is defined as the system of circles from $\mathcal C(p)$  that intersect $k^{2n}$. According to (\textbf{CS}), if the $n-$fold point of $k^{2n}$ coincide with $P_1$, the obtained surface $\mathcal {CS}(k^{2n}, p)$  is entirely spherical surface of order $2n$ having two real $n-$fold points in $P_1$ and $P_2$.

\begin{minipage}{7cm}
\begin{center}

\mbox{\includegraphics[width=7.5cm]{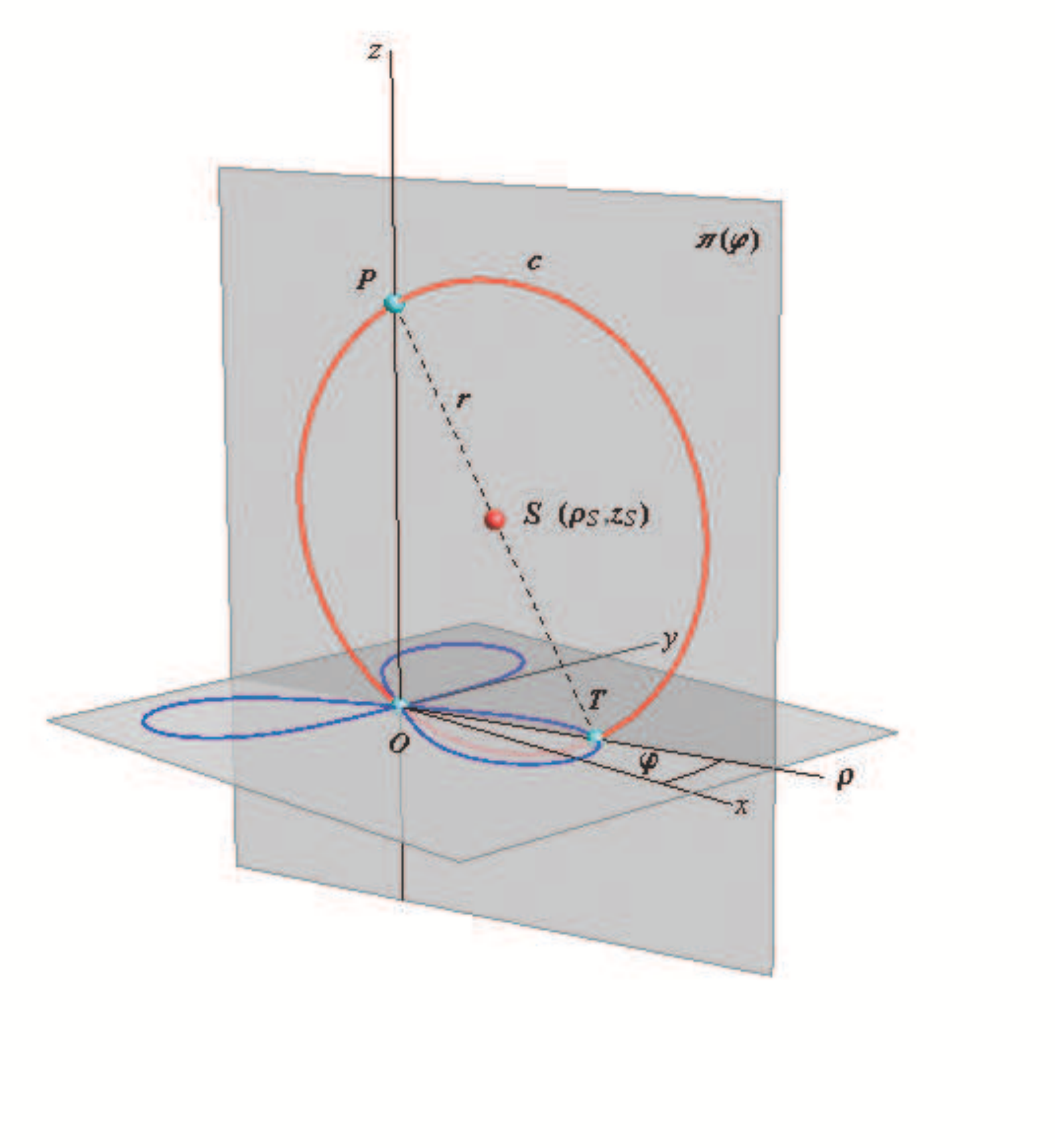}}\\

Figure 6

\end{center}
\end{minipage}
\begin{minipage}{7cm}
Here we offer a simpler construction of such surface.  
Let $k^{2n}$ be a curve in $xy-$plane having an $n-$fold point in the origin $O$, and let $P\neq O$ be a point on the axis $z$. In a plane $\pi(\varphi)$, such that   $z\subset\pi$ and $\varphi=\angle(\pi,x)$, a circle $c$ of radius  $PT$, where $T=\pi\cap k^{2n}$, is considered. See Figure~6.

If  $n$ is odd, there is a unique circle $c$ in every plane $\pi(\varphi)$. If $n$ is even, there could be two or none such circles.  For given $n$ and $p$ all circles $c$ determine observed surface  $\mathcal {CS}(k^{2n}, \frac{p}{2})$. 
\end{minipage}

\bigskip

In the plane $\pi(\varphi)$, with coordinates $(\rho,z)$, the circle  $c$ has an equation 
\begin{equation}\label{rho}
(\rho - \rho_S)^2+(z-z_S)^2=r^2,
\end{equation}
where $\rho_S=\frac{\rho_T}{2}$, $z_S=\frac{p}{2}$, $r^2=\frac{\rho_T^2}{4}+\frac{p^2}{4}$, and $\rho_T$ is given by (\ref{polar - curve1}).

Therefore, for $\varphi\in [0,2\pi)$, the equation (\ref{rho}) takes the form
\begin{equation}
\rho^2-\rho\sqrt[n]{\sin n\varphi}+z^2-pz=0
\end{equation}
which presents the equation of the surface $\mathcal {CS}(k^{2n}, \frac{p}{2})  $ in cylindrical  coordinates $(\rho,\varphi,z)$. 
If we raise it to the $n-$th power, the equation of the surface can be written as 

\begin{equation}\label{cilindricka2}
((\rho^2+z^2)+(-pz))^n =\rho^n\sin n\varphi.
\end{equation}

Since the correspondence between cylindrical  and Cartesian coordinates is given by  $\rho =\sqrt{x^2+y^2}$ and $\sin\varphi=\frac{y}{\sqrt{x^2+y^2}}$, if we use introduced notation $A_2=x^2+y^2+z^2$ and multiple-angle formulas \cite{math-multiple}

$$ \sin n\varphi= \left\{\begin{array}{ll}
 {(-1)^{(n-1)/2}\,\,T_n(\sin \varphi)}, & n \quad \textrm{odd},\\
{(-1)^{n/2-1}\cos\varphi\,\, U_{n-1}(\sin \varphi)}, & n \quad \textrm{even},
 \end{array}\right.
$$
where $T_n$ and $U_n$ are Chebyshev polynomials of the first and second kind, we can write equation  (\ref{cilindricka2}) in the following form:

\begin{equation}\label{implicit-ploha}
A_2^n+\sum_{j=1}^{n-1}\binom{n}{j}
(-pz)^j A_2^{n-j}=
\mathcal G^n(x,y)-(-pz)^n
\end{equation}

where

$$ \mathcal G^n(x,y)= \left\{\begin{array}{ll}
 {(-1)^{\frac{n-1}{2}}\left(\sqrt{x^2+y^2}\right)^n T_n\left(\frac{y}{\sqrt{x^2+y^2}}\right)}, & n \quad \textrm{odd},\\
{(-1)^{\frac{n}{2}-1}x\left(\sqrt{x^2+y^2}\right)^{n-1} U_{n-1}\left(\frac{y}{\sqrt{x^2+y^2}}\right)}, & n \quad \textrm{even}.
 \end{array}\right.
$$

The properties of the polynomials $T_n$ and $U_{n-1}$ (\cite{math-chebyshev1}, \cite{math-chebyshev2}), with  $x^2+y^2>0$, lead us to the conclusion that $\mathcal G^n$ are homogeneous polynomials  of degree $n$ in  $x$ and $y$. Therefore, the equation of constructed surface can be given in the following form:

\begin{equation}
A_2^n+\sum_{j=1}^{n-1}\binom{n}{j}
(-pz)^j A_2^{n-j}-
\mathcal T^n(x,y,z)=0,
\end{equation}
where $\mathcal T^n(x,y,z)=\mathcal G^n(x,y)-(-pz)^n$ is a homogeneous polynomial  of degree $n$ in  $x$, $y$ and $z$ determining the equation of tangent cone of the surface at its $n-$fold point $O$.  By using {\it Mathematica} functions {\it ChebyshevT} and {\it ChebyshevU} it is easy to obtain $T^n(x,y,z)$ for every $n$ and $p$.

The surface is symmetrical with respect to the plane  $z=\frac{p}{2}$. Thus, the tangent cone at the  $n-$fold point $P$ is symmetrical to the tangent cone at the point $O$ with respect to the same plane. Some examples of these surfaces are depicted in Figures 7 and 8.

\bigskip

\begin{center}
 \includegraphics[width=15cm]{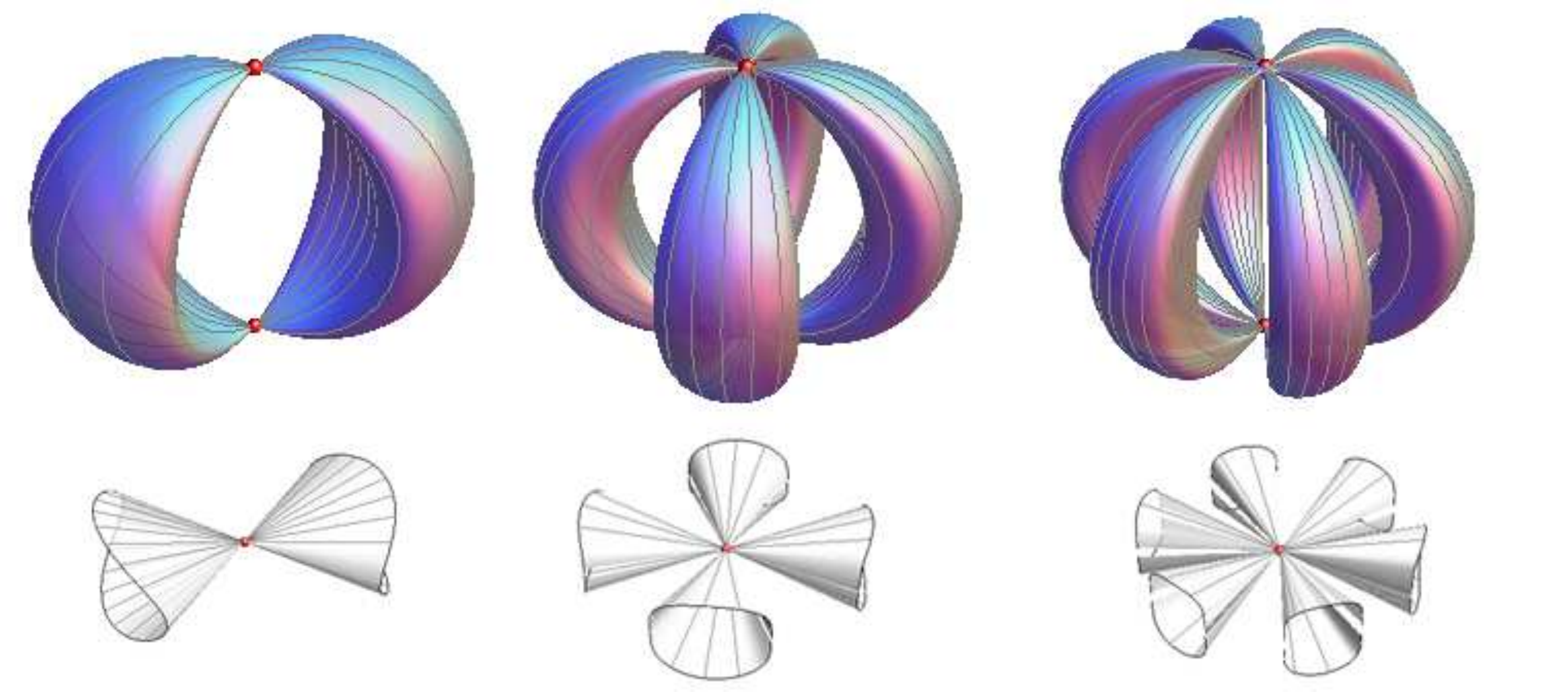}\\
\vspace{-0.5cm}
\begin{multicols}{3}
a

b

c

\end{multicols}
\end{center}
Figure 7: Three examples  of surfaces given by equation (\ref{implicit-ploha}) for $n=2,4,6$,  $p=2$ and their $n^{\mathrm{th}}$ degree tangent cones at the origin. The equations of tangent cones are: $x y-2 z^2=0$ (case a), $x^3 y-x y^3-4 z^4=0$ (case b) and $3 x^5 y-10 x^3 y^3+3 x y^5-32 z^6=0$ (case c).

\begin{center}
 \includegraphics[width=15cm]{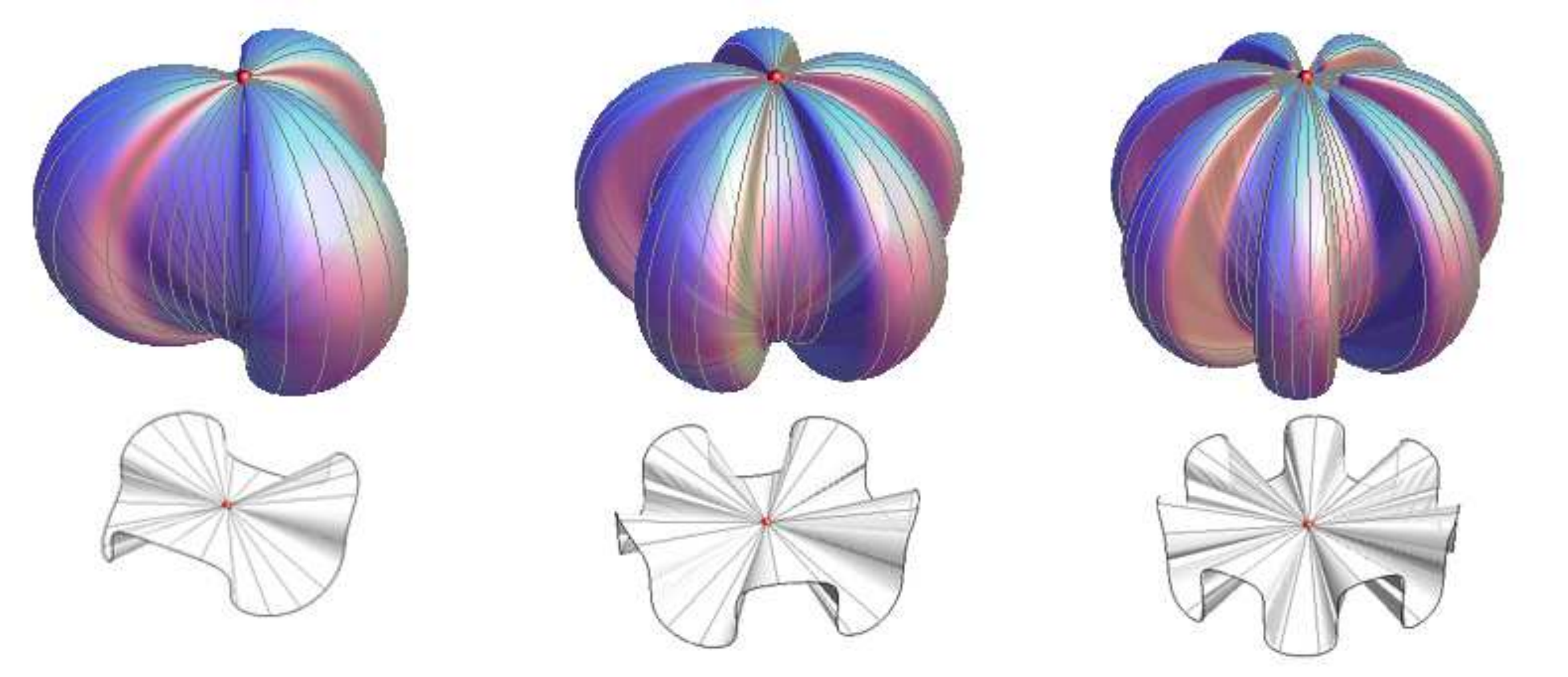}\\
\vspace{-0.5cm}
\begin{multicols}{3}
a

b

c

\end{multicols}
\end{center}
Figure 8: Three examples  of surfaces given by equation (\ref{implicit-ploha}) for $n=3,5,7$,  $p=2$ and their $n^{\mathrm{th}}$ degree tangent cones at the origin. The equations of tangent cones are: $x^2 y-y^3+8 z^3=0$ (case a), $x^4 y-10 x^2 y^3+y^5+32 z^5=0$ (case b) and $x^6 y-35 x^4 y^3+21 x^2 y^5-y^7+128 z^7=0$ (case c).

\end{document}